\documentstyle[amscd,amssymb,verbatim,12pt]{amsart}
\newcommand{\ch}{\mbox{\rm ch}}
\newcommand{\point}{\mbox{\rm pt.}}
\newcommand{\Def}{\mbox{\rm def}}
\newcommand{\asdim}{\mbox{\rm asdim}}
\newcommand{\Ricci}{\mbox{\rm Ricci}}

\newcommand{\Id}{\mbox{\rm Id}}

\newcommand{\Ker}{\mbox{\rm Ker}}

\newcommand{\Image}{\mbox{\rm Im}}

\newcommand{\Tr}{\mbox{\rm Tr}}

\newcommand{\area}{\mbox{\rm area}}
\newcommand{\ind}{\mbox{\rm ind}}
\newcommand{\rk}{\mbox{\rm rk}}

\newcommand{\HH}{\mbox{\rm H}}
\newcommand{\KK}{\mbox{\rm K}}
\newcommand{\uH}{\mbox{\rm H}_{(2)}}
\newcommand{\rH}{\overline{\mbox{\rm H}}_{(2)}}

\newcommand{\Dom}{\mbox{\rm Dom}}

\newcommand{\vol}{\mbox{\rm vol}}

\newcommand{\Isom}{\mbox{\rm Isom}}

\newcommand{\Int}{\mbox{\rm int}}

\newcommand{\R}{{\Bbb R}}
\newcommand{\C}{{\Bbb C}}

\newcommand{\Z}{{\Bbb Z}}
\theoremstyle{plain}
\newtheorem{definition}{Definition}

\newtheorem{lemma}{Lemma}
\newtheorem{proposition}{Proposition}
\newtheorem{corollary}{Corollary}

\newtheorem{conjecture}{Conjecture}
\numberwithin{equation}{section}
\renewcommand{\rm}{\normalshape}
\begin{document}
\title{The Zero-in-the-Spectrum Question}
\author{John Lott}
\address{Department of Mathematics\\
University of Michigan\\
Ann Arbor, MI  48109\\
USA}
\email{lott@@math.lsa.umich.edu}
\thanks{Research supported by NSF grant DMS-9403652}
\date{April 14, 1996}
\maketitle
\begin{abstract}
This is an expository article on the question of 
whether zero lies in the spectrum of
the Laplace-Beltrami operator acting on differential forms on a 
manifold.
\end{abstract}

\section{Introduction}
Let $M$ be a complete connected oriented Riemannian manifold.
The Laplace-Beltrami operator $\triangle_p$ acts on the
square-integrable $p$-forms on $M$. We asked the following question
in 1991 :\\ \\
{\bf Zero-in-the-Spectrum Question : }
Is zero always in the spectrum of
$\triangle_p$ for some $p$?\\

To our knowledge,
nobody has found a counterexample. The question was also raised by 
Gromov in the case of
a contractible manifold with a discrete cocompact group of isometries
\cite[p. 21]{Gromov (1993)}.

Being able to answer the above question is a first step toward
understanding the spectrum of the Laplace-Beltrami operator.
We would also like to be able to say
whether or not zero is in the spectrum of 
$\triangle_p$ for a given $p$. 
This problem is
partly topological in nature
and partly geometric, in a sense which will be made precise later.
In fact, it is equivalent to knowing the (unreduced) $L^2$-cohomology 
of $M$. The study of $L^2$-cohomology touches on many branches of mathematics,
including combinatorial group theory, topology, differential geometry and
algebraic geometry. It is most commonly considered
when $M$ is the universal cover of a compact manifold or when
$M$ is a finite-volume Hermitian locally symmetric space.  We refer to 
\cite{Lueck (1996),Pansu (1996)} and 
\cite{Stern (1993)} for surveys of these two cases.  In this
article we will instead emphasize general complete Riemannian manifolds
and give some partial positive answers to the 
zero-in-the-spectrum question.

The sections of the article are\\
1. Introduction\\
2. Definition of $L^2$-Cohomology\\
3. General Properties of $L^2$-Cohomology\\
4. Very Low Dimensions\\
4.1. One Dimension\\
4.2. Two Dimensions\\
5. Universal Covers\\
5.1. Big and Small Groups\\
5.2. Two and Three Dimensions\\ 
5.3. Four Dimensions\\
5.4. More Dimensions\\ 
6. Topologically Tame Manifolds\\

In what follows, all manifolds will be smooth, connected, oriented and of
positive dimension. All maps between manifolds
will be orientation-preserving. Unless otherwise indicated, all Riemannian
manifolds will be complete. 

We have tried to give as many complete proofs
as reasonably possible.
All unattributed results are of unknown origin or are due to the
author.  I thank Wolfgang L\"uck for conversations on some of the
topics discussed herein. I thank Marie-Claude Vergne for making the
figures. This article is based on lectures given
at the Troisi\`eme Cycle Romand held at Les Diablerets, Switzerland,
March, 1996.  I warmly thank Alain Valette and the other organizers
and participants of the meeting.
\section{Definition of $L^2$-Cohomology} \label{deff}
Let $M$ be as above. Let $\Lambda^p(M)$ denote the Hilbert space of
square-integrable $p$-forms on $M$. The completeness of $M$ enters in
one crucial way, in allowing us to integrate by parts on $M$ in the
sense of the following lemma.
\begin{lemma} (Gaffney \cite{Gaffney (1954)}) \label{parts}
Suppose that $\omega$, $\eta$, $d\omega$ and $d\eta$ are smooth
square-integrable differential forms on $M$. Then
\begin{equation}
\int_M d\omega \wedge \eta + (-1)^{deg(\omega)} \int_M \omega \wedge
d\eta = 0.
\end{equation}
\end{lemma}
\begin{pf}
We claim that there is a sequence
$\{\phi_i\}_{i=1}^\infty$ of compactly-supported functions on $M$ with 
the properties that\\
1. There is a constant $C > 0$ such that for all $i$ and almost all $m \in M$,
$|\phi_i(m)| \le C$ and $|d\phi_i(m)| \le C$.\\
2. For almost all $m \in M$, $\lim_{i \rightarrow \infty} \phi_i(m) = 1$ and
$\lim_{i \rightarrow \infty} |d\phi_i(m)| = 0$.

To construct the sequence $\{\phi_i\}_{i=1}^\infty$, let $m_0$ be a
basepoint in $M$. Let $f \in C^\infty_0([0, \infty))$ be a nonincreasing
function such that if $x \in [0,\frac12]$ then $f(x) = 1$. Put
$\phi_i(m) = f \left(\frac1i d(m_0,m)\right)$. 
This gives the desired sequence. The completeness of $M$ ensures that
$\phi_i$ is compactly-supported.
Note that $\phi_i$ is {\it a priori} only a Lipschitz function, but this is
good enough for our purposes.

Using Lebesgue Dominated Convergence and the fact that we can integrate
by parts for compactly-supported forms, we have
\begin{align}
\int_M d\omega \wedge \eta + (-1)^{deg(\omega)} \int_M \omega \wedge
d\eta & = \int_M d(\omega \wedge \eta) = \lim_{i \rightarrow \infty}
\int_M \phi_i \: d(\omega \wedge \eta) \\
& = - \lim_{i \rightarrow \infty}
\int_M d\phi_i \wedge \omega \wedge \eta = 0. \notag
\end{align}
This proves the lemma.
\end{pf}  

Let $d^*$ be the formal adjoint to $d$. Using Lemma \ref{parts},
one can construct a self-adjoint operator $\triangle = d d^* + d^* d$ 
acting on $\Lambda^*(M)$, with domain
\begin{equation}
\Dom(\triangle) = \{ \omega \in \Lambda^*(M) : d\omega, \: d^* \omega, \:
d d^* \omega \: \text{ and } d^* d \omega \text{ are square-integrable}\}.
\end{equation} 
Let $\triangle_p$ denote the restriction of $\triangle$ to $\Lambda^p(M)$.
The spectrum $\sigma(\triangle_p)$ of $\triangle_p$
is a closed subset of $[0, \infty)$.
\begin{lemma} \label{kernel}
The kernel of $\triangle_p$ is
$ \{ \omega \in \Lambda^p(M) : d \omega = d^* \omega = 0\}$.
\end{lemma}
\begin{pf}
Clearly $\{ \omega \in \Lambda^p(M) : d \omega = d^* \omega = 0\} \subseteq
\Ker(\triangle_p)$. If $\omega \in \Ker(\triangle_p)$ then by elliptic
regularity, $\omega$ is smooth. Using integration
by parts,
$0 = \langle \omega, \triangle_p \omega \rangle = \langle d \omega, d \omega
\rangle + \langle d^* \omega, d^* \omega \rangle$, so
$d\omega = d^* \omega = 0$.
\end{pf}
{\bf Warning : } Unlike what happens with compact manifolds, it is 
possible that $\Ker(\triangle_p) = 0$ but nevertheless
$0 \in \sigma(\triangle_p)$. The simplest example of this is when
$M = \R$ and $p = 0$. By Lemma \ref{kernel}, $\Ker(\triangle_0)$ consists of
square-integrable functions $f$ on $\R$ such that $df = 0$. Clearly
the only such function is the zero function. 
However, under Fourier transform, $\triangle_0$ is
equivalent to the multiplication operator by $k^2$ on $L^2(\R)$ and
hence $\sigma(\triangle_0) = [0, \infty)$. \\ \\
{\bf Examples : } We now give $\sigma(\triangle_p)$ for simply-connected
space forms.\\ \\
1. $M$ is the standard sphere $S^n$. From \cite{Gallot-Meyer (1975)},
\begin{equation}
\sigma(\triangle_p) = \{ (k+p) (k+n+1-p)\}_{k=0}^\infty \cup
\{ (k+p+1) (k+n-p)\}_{k=0}^\infty.
\end{equation}
(See Fig. 1.)
The details of the spectrum are not important for us.  We only wish to 
note that $\sigma(\triangle_p)$ is discrete, and 
$0 \in \sigma(\triangle_p)$ if $p = 0$ or $p=n$. These
statements are a consequence of the fact that $M$ is closed.  Namely,
if $M^n$ is any closed 
Riemannian manifold then $\sigma(\triangle_p)$
is discrete and $\Ker(\triangle_p) \cong \HH^p(M; \C)$. In particular,
$\Ker(\triangle_0) \cong \HH^0(M; \C) = \C$ consists of the constant
functions and $\Ker(\triangle_n) \cong \HH^n(M; \C) = \C$ consists of
multiples of the volume form.\\ \\
2. $M$ is the standard Euclidean space $\R^n$. 
As the $p$-forms on $\R^n$ consist of $\binom{n}{p}$ copies
of the functions, it is enough to consider $\sigma(\triangle_0)$. By
Fourier analysis, $\sigma(\triangle_0) = [0, \infty)$. Thus
$\sigma(\triangle_p) = [0, \infty)$ for all $0 \le p \le n$. 
(See Fig. 2.) Note that
$\Ker(\triangle_p) = 0$ for all $p$.\\ \\
3. $M$ is the hyperbolic space $H^{2n}$. From \cite{Donnelly (1981)},
\begin{equation}
\sigma(\triangle_p) = 
\begin{cases}
\left[ \frac{(2n-2p-1)^2}{4}, \infty \right) & \text{ if } 0 \le p \le n-1,\\
\{0\} \cup \left[ \frac{1}{4}, \infty \right) & \text{ if } p=n, \notag \\
\left[ \frac{(2p-2n-1)^2}{4}, \infty \right) &
\text{ if } n+1 \le p \le 2n. \notag
\end{cases}
\end{equation}
(See Fig. 3.)
There is an infinite-dimensional kernel to $\triangle_n$. Otherwise, the
spectrum is strictly bounded away from zero.\\ \\
4. $M$ is the hyperbolic space $H^{2n+1}$. From \cite{Donnelly (1981)},
\begin{equation}
\sigma(\triangle_p) = 
\begin{cases}
\left[ \frac{(2n-2p)^2}{4}, \infty \right) & \text{ if } 0 \le p \le n,\\ 
\notag \\
\left[ \frac{(2p-2n-2)^2}{4}, \infty \right) & 
\text{ if } n+1 \le p \le 2n+1. \notag
\end{cases}
\end{equation}
(See Fig. 4.)
For all $p$, $\Ker(\triangle_p) = 0$. The continuous spectrum extends down
to zero in degrees $n$ and $n+1$, and is strictly bounded away from zero
in other degrees.\\
{\bf End of Examples}\\

Comparing Figures 1-4, 
the spectra do not have much in common.
However, one common feature is that zero lies in $\sigma(\triangle_p)$ for some
$p$, although for different reasons
in the different cases.  In Figure 1, it is because
$\triangle_0$ has a nonzero finite-dimensional kernel.  In Figure 2,
it is because zero lies in the continuous spectrum of $\triangle_p$ for
all $p$.  In Figure 3, it is because $\triangle_n$ has an infinite-dimensional
kernel.  And in Figure 4, it is because zero lies in the continuous
spectrum of $\triangle_p$ for $p = n$ and $p = n+1$.

The above examples, along with others, 
motivate the zero-in-the-spectrum question.  One
can pose the question for various classes of manifolds, such as\\
1. Complete Riemannian manifolds.\\
2. Complete Riemannian manifolds of bounded geometry, meaning that
the injectivity radius is positive and the sectional curvature $K$
satisfies $|K| \le 1$.\\
3. Uniformly contractible Riemannian manifolds, meaning that for all
$r > 0$, there is an $R(r) \ge r$ such that for all $m \in M$, the
metric ball $B_r(m)$ can be contracted to a point within $B_{R(r)}(m)$.\\
4. Universal covers of closed Riemannian manifolds.\\
5. Universal covers of closed aspherical Riemannian manifolds.\\

There are obvious inclusions
$\begin{array}{lllll}
5 & \subset & 4 & \subset & 2\\
\cap  &         & & & \cap \notag \\
3 & & \subset & & 1. \notag
\end{array}$
As we shall discuss, there are some reasons to believe that the answer to the 
zero-in-the-spectrum question is ``yes'' in class 5, but the evidence 
for a ``yes'' answer in class 1 consists mainly of a lack of counterexamples.

In order to make the study of the spectrum of $\triangle_p$ more precise,
the Hodge decomposition
\begin{equation} \label{Hodge}
\Lambda^p(M) = \Ker(\triangle_p) \oplus \overline{\Image(d)} \oplus
\Lambda^p(M)/\Ker(d)
\end{equation}
is useful.
The operator $\triangle_p$ decomposes with respect to
(\ref{Hodge}) as a direct sum of three operators.
If we know the spectrum of the Laplace-Beltrami operator on all forms
of degree less than $p$ then the new information in degree $p$ consists
of $\Ker(\triangle_p)$ and the spectrum of $\triangle_p$ on 
$\Lambda^p(M)/\Ker(d)$. So we can ask the more precise questions :\\
1. What is $\dim(\Ker(\triangle_p))$?\\
2. Is zero in $\sigma \left( \triangle_p \text{ on } 
\Lambda^p(M)/\Ker(d) \right)$? \\

By its definition, $\triangle_p$ involves the first derivatives of the
metric tensor.  We now show that the answer to the zero-in-the-spectrum
question only depends on the $C^0$-properties of the metric tensor.
To do so, we reformulate the question in terms of $L^2$-cohomology.
Define a subspace $\Omega^p(M)$ of $\Lambda^p(M)$ by
\begin{equation}
\Omega^p(M) = \{ \omega \in \Lambda^p(M) : d\omega 
\text{ is square-integrable}\},
\end{equation}
where $d\omega$ is initially interpreted in a distributional sense.
The subspace $\Omega^p(M)$ is cooked up so that we have a cochain complex
\begin{equation} \label{cochaincomplex}
\ldots \stackrel{d_{p-1}}{\longrightarrow} \Omega^p(M)
\stackrel{d_{p}}{\longrightarrow} \Omega^{p+1}(M)
\stackrel{d_{p+1}}{\longrightarrow} \ldots
\end{equation}
\begin{lemma}
$\Ker(d_p)$ is a closed subspace of $\Lambda^p(M)$.
\end{lemma}
\begin{pf}
Suppose that $\{ \eta_i \}_{i=1}^\infty$ is a sequence in $\Ker(d_p)$ 
which converges to $\omega \in \Lambda^p(M)$ in an $L^2$-sense.
We must show that the distributional form $d \omega$ vanishes.
Given a smooth compactly-supported $(p + 1)$-form $\rho$,
we have 
\begin{equation}
\langle d \omega, \rho \rangle = \langle \omega, d^* \rho \rangle =
\lim_{i \rightarrow \infty} \langle \eta_i, d^* \rho \rangle =
\lim_{i \rightarrow \infty} \langle d \eta_i, \rho \rangle = 0.
\end{equation}
The lemma follows.
\end{pf} 
\begin{definition} \label{L^2coho}
The $p$-th unreduced $L^2$-cohomology group of $M$ is
$\uH^p(M) = \Ker(d_p)/\Image(d_{p-1})$. The $p$-th reduced $L^2$-cohomology
group of $M$ is 
$\rH^p(M) = \Ker(d_p)/\overline{\Image(d_{p-1})}$, a Hilbert space.
\end{definition}

The square-integrability condition on the forms should be thought of as
a global decay condition, not as a local regularity condition. One
can also compute $\uH^*(M)$ using a complex as in (\ref{cochaincomplex}) 
where the forms are additionally required to be smooth 
\cite[Prop. 9]{Lott (1996)}.

There is an obvious surjection $i_p : \uH^p(M) \rightarrow \rH^p(M)$. 
Clearly $i_p$ is an isomorphism if and only if $d_{p-1}$ has closed image.

\begin{proposition}
1. $\Ker(\triangle_p) \cong \rH^p(M)$.\\
2. $0 \notin \sigma \left( \triangle_p \text{ on } 
\Lambda^p(M)/\Ker(d) \right)$ if and only if $i_{p+1}$ is an isomorphism.
\end{proposition} 
\begin{pf}
1. Using Lemma \ref{kernel}, we have 
\begin{equation}
\Ker(\triangle_p) = \{ \omega \in \Lambda^p(M) : d \omega = d^* \omega = 0\}
= \Ker(d_p) \cap \overline{\Image(d_{p-1})}^\perp \cong \rH^p(M).
\end{equation}
The first part of the proposition follows.\\
2. Suppose first that 
$\triangle_p$ has a bounded inverse on $\Lambda^p(M)/\Ker(d)$.
Given $\mu \in \Lambda^p(M)$, let
$\overline{\mu}$ denote its class in $\Lambda^p(M)/\Ker(d)$. 
Define an operator $S$ on smooth compactly-supported $(p+1)$-forms by
$S (\omega) = d \triangle_p^{-1} 
\overline{d^* \omega}$.
Then $S$ extends to a bounded operator on $\Lambda^{p+1}(M)$. Let
$\{\eta_i\}_{i=1}^\infty$ be a sequence in $\Omega^p(M)$ such that
$\lim_{i \rightarrow \infty} d \eta_i = \omega$ for some
$\omega \in \Lambda^{p+1}(M)$. 
Then for each $i$, we have $d \eta_i = S(d \eta_i)$ and so 
$\omega = S(\omega)$. Thus $\omega \in \Image(d)$ and so
$\Image(d)$ is closed.

Now suppose that 
$\triangle_p$ does not have a bounded inverse on $\Lambda^p(M)/\Ker(d)$.
Then there is a sequence of positive numbers
$r_1 > s_1 > r_2 > s_2 > \ldots$ tending towards zero and
an orthonormal sequence $\{\eta_i\}_{i=1}^\infty$ in
$\Lambda^p(M)/\Ker(d)$ such that with respect to the spectral projection $P$
of $\triangle_p$ (acting on $\Lambda^p(M)/\Ker(d)$), 
$\eta_i \in \Image(P([s_i, r_i]))$.
Put $\lambda_i = \| d \eta_i \|$.
Then $\lim_{i \rightarrow \infty} \lambda_i = 0$.
Let $\{ c_i \}_{i=1}^\infty$ be a sequence in $\R^+$ such that
$\sum_{i = 1}^\infty c_i^2 = \infty$ and
$\sum_{i = 1}^\infty c_i \lambda_i < \infty$.  Put $\omega =
\sum_{i=1}^\infty c_i d \eta_i$. Then $\omega \in 
\overline{\Image(d)}$.  Suppose that $\omega = d \mu$ for some
$\mu \in \Omega^p(M)$.
By the spectral theorem, we must have $\overline{\mu} = 
\sum_{i = 1}^\infty c_i \eta_i$.
However, this is not square-integrable.
Thus $\Image(d)$ is not closed. The proposition follows.
\end{pf}

\begin{corollary}
Zero does not lie in $\sigma(\triangle_p)$ for any $p$ if and only if
$\uH^p(M)$ = 0 for all $p$, i.e. if the complex
(\ref{cochaincomplex}) is contractible.
\end{corollary}

So a counterexample to the zero-in-the-spectrum question would consist
of a manifold $M$ whose complex (\ref{cochaincomplex}) is contractible.
By way of comparison, recall that the compactly-supported 
complex-valued cohomology of $M$
 is computed by a cochain complex similar to (\ref{cochaincomplex}), except 
using  compactly-supported smooth forms.  As $\HH_c^{dim(M)}(M; \C) \neq 0$,
this latter complex is never contractible.  
And the ordinary complex-valued cohomology of $M$
is computed by a cochain complex similar to (\ref{cochaincomplex}), except 
using smooth forms without any decay conditions. Again, as
$\HH^0(M; \C) \neq 0$, this latter complex is never contractible.

If $M$ is closed then $\rH^*(M)$ is independent of the Riemannian metric
on $M$. This is no longer true if $M$ is not closed - consider $\R^2$ and
$H^2$. However, the $L^2$-cohomology groups of $M$ do have some
invariance properties which we now discuss.

\begin{definition}
Riemannian manifolds $M$ and $M^\prime$ 
are biLipschitz diffeomorphic if
there is a diffeomorphism $F : M \rightarrow M^\prime$ and a constant $K > 0$
such that the Riemannian metrics $g$ and $g^\prime$ satisfy the pointwise 
inequality
\begin{equation}
K^{-1} g \le F^* g^\prime \le K g.
\end{equation}
\end{definition}
If $M$ and $M^\prime$ are biLipschitz diffeomorphic then their reduced
and unreduced $L^2$-cohomology groups are isomorphic, as the Riemannian
metric only enters in the complex (\ref{cochaincomplex}) in determining
which forms are square-integrable. Thus the answer to the 
zero-in-the-spectrum question only depends on the biLipschitz diffeomorphism
class of $M$. More generally, we can consider a category whose objects are
Lipschitz Riemannian manifolds and 
whose morphisms are Lipschitz
maps.  Then the reduced and unreduced $L^2$-cohomology groups are
Lipschitz-homotopy-invariants.

Note that $L^2$-cohomology groups are not coarse quasi-isometry invariants.
For example, any closed manifold is coarsely quasi-isometric to a point,
but its $L^2$-cohomology is the same as its ordinary complex-valued 
cohomology,
which may not be that of a point. However, some aspects of $L^2$-cohomology
only depend on the large-scale geometry of the manifold.
\begin{proposition} \cite[Prop. 12]{Lott (1996)} \label{outside}
If $M$ and $M^\prime$ are isometric
outside of compact sets then\\
1. $\Ker(\triangle_p)$ is finite-dimensional on $M$ if and only if it is
finite-dimensional on $M^\prime$.\\
2. Zero is in $\sigma \left( \triangle_p \text{ on } \Lambda^p/\Ker(d) \right)$
on $M$ if and only if the same statement is true on $M^\prime$.
\end{proposition}

Consider uniformly contractible Riemannian manifolds of bounded geometry.
If two such manifolds are coarsely quasi-isometric then they are
Lipschitz-homotopy-equivalent and hence their $L^2$-cohomology groups
are isomorphic \cite[p. 219]{Gromov (1993)}. The next proposition gives
an extension of this result in which uniform contractibility is replaced by
uniform vanishing of cohomology, the latter being defined as follows.

\begin{definition}
We say that $\HH^j(M; \C)$ vanishes uniformly if
for all $r > 0$, there is an $R(r) \ge r$ such that for all $m \in M$,
\begin{equation}
\Image \left(H^j(B_{R(r)}(m); \C) \rightarrow H^j(B_{r}(m); \C) \right) = 0.
\end{equation}
\end{definition}

\begin{proposition} \label{pansu} (Pansu \cite{Pansu (1995)})
Consider a
Riemannian manifold $M$ of bounded geometry such that for some $k > 0$,
$\HH^j(M; \C)$ vanishes uniformly for $1 \le j \le k$.
Then within the class of such manifolds,\\
1. $\rH^p(M)$ and $\uH^p(M)$ 
are coarse quasi-isometry invariants for
$0 \le p \le k$.\\
2. $\Ker(\rH^{k+1}(M) \rightarrow \HH^{k+1}(M; \C))$ and
$\Ker(\uH^{k+1}(M) \rightarrow \HH^{k+1}(M; \C))$ 
are coarse quasi-isometry invariants.
\end{proposition}
\section{General Properties of $L^2$-Cohomology}

In this section we give some general results about the $L^2$-cohomology
of complete Riemannian manifolds.  First, we give a useful sufficient
condition for the reduced $L^2$-cohomology to be nonzero.
\begin{proposition} \label{compactsup}
For all $p$, $\Image \left(\HH^p_c(M; \C) \rightarrow \HH^p(M; \C) \right)$ 
injects into $\rH^p(M)$.
\end{proposition}
\begin{pf}
Suppose that $\omega$ is a smooth compactly-supported closed $p$-form
which represents a nonzero class in $\HH^p(M; \C)$.
By Poincar\'e duality, there is a smooth compactly-supported closed 
$(\dim(M) - p)$-form $\rho$ such that $\int_M \omega \wedge \rho \neq 0$. 

As $\omega$ is compactly-supported, it is square-integrable and so represents
an element $[\omega]$ of $\rH^p(M)$. Suppose that $[\omega] = 0$. Then
there is a sequence $\{\eta_i\}_{i=1}^\infty$ in $\Omega^{p-1}(M)$ such
that $\omega = \lim_{i \rightarrow \infty} d \eta_i$, where the limit
is in an $L^2$-sense.  It follows that
\begin{equation}
\int_M \omega \wedge \rho = 
\lim_{i \rightarrow \infty} \int_M d\eta_i \wedge \rho =
\lim_{i \rightarrow \infty} \int_M d(\eta_i \wedge \rho) = 0,
\end{equation}
which is a contradiction.  Thus $[\omega] \neq 0$.
\end{pf}
\begin{corollary} \label{cor}
Let $N^{4k}$ be a compact manifold-with-boundary 
with nonzero signature.  Then if $M$ is any complete Riemannian manifold
which is diffeomorphic to $\Int(N)$, $\rH^{2k}(M) \neq 0$.
\end{corollary}
\begin{pf}
By definition, the signature of $N$ is the signature of the intersection
form on  
\begin{equation}
\Image \left(\HH^{2k}(N, \partial N; \C) \rightarrow 
\HH^{2k}(N; \C) \right) \cong
\Image \left(\HH^{2k}_c(M; \C) \rightarrow \HH^{2k}(M; \C) \right).
\end{equation}
If the signature of $N$ is nonzero then 
$\Image \left(\HH^{2k}_c(M; \C) \rightarrow \HH^{2k}(M; \C) \right)$ must
be nonzero. The corollary follows from Proposition \ref{compactsup}.
\end{pf}
{\bf Example : } Let $N$ be $\C P^2$ with a small $4$-ball removed. Then
$N$ satisfies the hypothesis of Corollary \ref{cor}.\\

We now show that the middle-dimensional reduced $L^2$-cohomology is a
conformal invariant of $M$.

\begin{proposition} \label{conformal}
If $M^{2k}$ is even-dimensional then $\Ker(\triangle_k)$ 
is conformally-invariant.
\end{proposition}
\begin{pf}
Suppose that $g$ and $e^\phi g$ are conformally equivalent Riemannian
metrics on $M$, with $\phi \in C^\infty(M)$.
We use the fact that the action of the Hodge duality operator
$*$ on $\Lambda^k(M)$ is independent of
$\phi$. If $\omega$ is a $k$-form on $M$, its $L^2$-norm $\int_M \omega
\wedge *\omega$ is independent of $\phi$. Thus the Hilbert space $\Lambda^k(M)$
is independent of $\phi$. Furthermore,
\begin{align}
\Ker(\triangle_k) & = \{ \omega \in \Lambda^k(M) : d\omega = d^* \omega = 0\}\\
& = \{ \omega \in \Lambda^k(M) : d\omega = \pm *d* (\omega) = 0\} \notag \\
& = \{ \omega \in \Lambda^k(M) : d\omega = d * (\omega) = 0\}
\end{align}
is independent of $\phi$.
\end{pf}
\noindent
{\bf Example : } Take $M = H^2$. Then $M$ is conformally equivalent to
a Euclidean disk $D$. The harmonic square-integrable $1$-forms on $D$ are
of the form $f_1(x,y) dx + f_2(x,y) dy$, where $f_1$ and $f_2$ are
square-integrable harmonic functions on $D$. There is clearly an infinite
number of such functions, and so $\dim(\rH^1(H^2)) = \infty$.
The same argument applies to $H^{2k}$, to give
$\dim(\rH^k(H^{2k})) = \infty$.\\  

In the case of functions, one has a good control of
when zero is in the spectrum of the Laplacian.
\begin{lemma}
$\Ker(\triangle_0) \neq 0$ if and only if $\vol(M) < \infty$.
\end{lemma}
\begin{pf}
If $\vol(M) < \infty$ then the constant functions on $M$ are
square-integrable and harmonic. Conversely, if $f \in \Ker(\triangle_0)$
then by Lemma \ref{kernel}, $f$ is constant. If $f$
is nonzero and square-integrable then $\vol(M) < \infty$.
\end{pf}

\begin{definition} \label{open} $M$ is open at infinity if
there is a constant $C > 0$ such that for all domains $D$ in $M$
with smooth compact closure, $\frac{\area(\partial D)}{\vol(D)} \ge C$.
\end{definition}
\noindent
{\bf Examples : }\\
1. $\R^n$ is not open at infinity, as can be seen by taking large balls
for $D$.\\
2. $H^n$ is open at infinity.
\begin{proposition} (Buser \cite{Buser (1982)}) \label{Buser}
Let $M$ have infinite volume.
Suppose that there is a constant
$c \ge 0$ such that $\Ricci_M \ge - c^2$.  Then $0 \notin \sigma(\triangle_0)$
if and only if $M$ is open at infinity.
\end{proposition}
\begin{pf}
1. Suppose that $M$ is open at infinity. By Cheeger's inequality,
\begin{equation}
\inf(\sigma(\triangle_0)) \ge \inf_D \frac14 \left(
\frac{\area(\partial D)}{\vol(D)} \right)^2 > 0.
\end{equation}
2. Suppose that $M$ is not open at infinity.  The bottom of the
spectrum of $\triangle_0$ is given in terms of Rayleigh quotients by
\begin{equation}
\inf(\sigma(\triangle_0)) = \inf_{f \ne 0} 
\frac{\int_M |df|^2}{\int_M f^2},
\end{equation}
where $f$ ranges over compactly-supported Lipschitz functions on $M$.
We want to find compactly-supported Lipschitz functions on $M$ of
arbitrarily small Rayleigh quotient. By assumption, for all $\epsilon > 0$
there is a domain $D$ such that 
$\frac{\area(\partial D)}{\vol(D)} \le \epsilon$.
Put 
\begin{equation}
N_1(\partial D) = \{ m \in M : m \notin D \text{ and }
d(m, \partial D) \le 1\}.
\end{equation}
Define a function $f$, which
approximates the characteristic function of $D$, by
\begin{equation}
f(m) =
\begin{cases}
1 & \text{ if } m \in D \\
1 - d(m, \partial D) & \text{ if } m \in N_1(\partial D) \\
0 & \text{ if } m \notin D$ and $m \notin N_1(\partial D).
\end{cases}
\end{equation}
Clearly $\int_M f^2 \ge \vol(D)$. As $f$ has nonzero gradient
only in $N_1(\partial D)$, where
$|df| = 1$ almost everywhere, 
we have $\int_M |df|^2 = \vol(N_1(\partial D))$. 
If $D$ is nice and round then we expect that 
\begin{equation} \label{spiky}
\vol(N_1(\partial D)) \sim \area(\partial D)
\end{equation}
and the Rayleigh quotient $\frac{\int_M |df|^2}{\int_M f^2}$ will be
comparable to $\epsilon$.

The only problem with this argument is that $D$ may not be nice and round,
but may have long thin legs coming out of it, like an octupus.  
Then (\ref{spiky}) may not be valid. The content of \cite{Buser (1982)} is
that if this is the case, we can cut the legs off of $D$ to
come up with a new domain for which the above heuristic argument is valid.
It is in this step that the lower bound on the Ricci curvature is used.
We refer to \cite{Buser (1982)} for details.
\end{pf}

\begin{corollary} (Brooks \cite{Brooks (1981)}) \label{Brooks}
Let $M$ be a normal covering of a compact manifold $X$ with covering
group $\Gamma$. Then $0 \in \sigma(\triangle_0)$ on $M$ if and only if
$\Gamma$ is amenable.
\end{corollary}
\begin{pf} If $\Gamma$ is finite then $0 \in \sigma(\triangle_0)$ and
$\Gamma$ is amenable.  If $\Gamma$ is infinite then by Proposition
\ref{Buser}, $0 \in \sigma(\triangle_0)$ if and only if $M$ is not open at
infinity.  Let $S$ be a finite set of generators of $\Gamma$. Let $G$ be
the Cayley graph of $\Gamma$, constructed using $S$. There is a notion of
$G$ being open at infinity which is similar to Definition \ref{open}.
As $M$ is
coarsely quasi-isometric to $G$, $M$ is not open at infinity if and only if
$G$ is not open at infinity.  However, this is one of
the characterizations of amenability of $\Gamma$.
\end{pf}

We now prove a result about manifolds which, roughly
speaking, are at least as large as Euclidean space.
\begin{definition}
$M$ is hyperEuclidean if there is a proper distance-nonincreasing map
$F : M \rightarrow \R^{\dim(M)}$ of nonzero degree.
\end{definition}
\noindent
{\bf Remarks : } \\
1. A map is proper if preimages of compact sets are compact.
Instead of requiring that $F$ be distance-nonincreasing, we could
require that $F$ have a finite Lipschitz constant.  
By postcomposing $F$ with a 
dilatation of $\R^{\dim(M)}$, the two conditions are equivalent.\\
2. If $M$ is hyperEuclidean then a compactly-supported modification of
$M$ is also hyperEuclidean.\\
3. Examples of hyperEuclidean manifolds are given by simply-connected
nonpositively-curved manifolds $M$. Namely, fix $m_0 \in M$ and
put $F = \exp_{m_0}^{-1}$.\\
4. There was once a conjecture that all uniformly contractible manifolds
are hyperEuclidean (with a degree-one map to $\R^{\dim(M)}$), 
but this turns out to be wrong 
\cite{Dranishnikov-Ferry-Weinberger (1994)}. There is still an open conjecture
that a uniformly contractible manifold of bounded geometry is hyperEuclidean,
and in particular, that the universal cover of an aspherical closed 
manifold is hyperEuclidean.  
\begin{proposition} (Gromov \cite[p. 238]{Gromov (1993)}) \label{hyper}
If $M$ is hyperEuclidean then $0 \in \sigma(\triangle_p)$ for some $p$.
\end{proposition}
\begin{pf}
Put $n = \dim(M)$.
First, suppose that $n$ is even.  We will construct a vector bundle $E$
with connection on $\R^n$ which is topologically nontrivial but
analytically trivial, in a sense which will be made precise.
Then assuming that zero is not in the spectrum of $M$, we will
apply the relative index theorem to $F^* E$ in order to get a
contradiction.

Recall that $\KK^0(S^n) = \Z \oplus \Z$. If ${\cal E}$ is a (virtual) 
vector bundle
on $S^n$, the two $\Z$ factors correspond to $\rk({\cal E})$ and 
$\int_{S^n} \ch({\cal E})$, respectively. This means that for some 
$N > 0$, there is a complex $\C^N$-bundle ${\cal E}$ on $S^n$ with 
$\int_{S^n} \ch({\cal E}) \ne 0$. Fixing a point $\infty \in S^n$,
we can trivialize ${\cal E}$ in a neighborhood of $\infty$. Furthermore,
we can put a Hermitian metric and Hermitian connection on ${\cal E}$
so that the connection is flat in a neighborhood of $\infty$.

Let $E$ be the restriction of ${\cal E}$ to $\R^n = S^n - \{\infty\}$.
Let $\nabla$ be the restriction of the Hermitian connection on ${\cal E}$ to
$\R^n$.
Then $E$ is trivialized outside of a compact set 
$K \subset \R^n$ and $\nabla$ is flat outside of $K$.

As $\R^n$ is contractible, there is an isomorphism of Hermitian vector bundles
$i : \R^n \times \C^N \rightarrow E$. Then $i^* \nabla$
can be considered to be a
$u(N)$-valued $1$-form $\omega$ on $\R^n$. The curvature of $\omega$ is the
$u(N)$-valued $2$-form 
$\Omega = d \omega + \omega^2$. The nontriviality of ${\cal E}$ translates to
the facts that\\
1. $\Omega$ vanishes outside of $K$ and\\
2. The de Rham cohomology class of the compactly-supported form 
$\Tr \left( e^{- \frac{\Omega}{2\pi i}} \right) - N$ is a nonzero multiple
of the fundamental class $[\R^n] \in \HH^n_c(\R^n; \R)$.\\

In fact, we can take $\omega$ to have a finite $L^\infty$-norm
$\parallel \omega \parallel_\infty$. For
example, if $n=2$, take $N = 1$. 
Let $f \in C^\infty_0([0, \infty))$ be a nonincreasing
function such that if $x \in [0,\frac12]$ then $f(x) = 1$. Put 
$\omega = - i (1-f(r)) d\theta$. Then 
\begin{equation}
\Omega = d \omega = i f^\prime(r) dr \wedge d\theta.
\end{equation} 
We have
$\parallel \omega \parallel_\infty \: = \sup_{r \ge 0} \frac{1-f(r)}{r}$ and 
$\int_{\R^2} \left[ \Tr \left( e^{- \frac{\Omega}{2\pi i}} \right) - 1
\right] = 1$.

Returning to the case of general even $n$, for $\epsilon > 0$, let
$\Phi_\epsilon : \R^n \rightarrow \R^n$ be the map $\Phi_\epsilon({\bf x}) =
\epsilon {\bf x}$. Put $\omega_\epsilon = \Phi_\epsilon^* \omega$ and
$\Omega_\epsilon = d \omega_\epsilon + \omega_\epsilon^2$. Then
\begin{equation}
\parallel \omega_\epsilon \parallel_\infty \: = \epsilon
\parallel \omega \parallel_\infty \text{ and } 
\int_{\R^n} \left[ \Tr \left( e^{- \frac{\Omega_\epsilon}{2\pi i}} \right) 
-N \right] =
\int_{\R^n} \left[ \Tr \left( e^{- \frac{\Omega}{2\pi i}} \right) - N
\right] \ne 0.
\end{equation}

We now turn our attention to $M$. Suppose that 
$0 \notin \sigma(\triangle_p)$ for all $p$. Consider the self-adjoint
operator $d + d^*$ on 
$\Lambda^*(M)$. As $(d + d^*)^2 = \triangle$, it follows that 
$0 \notin \sigma(d + d^*)$. In other words, $d + d^*$ is $L^2$-invertible.
Define an operator $\mu$ on $\Lambda^*(M)$ by saying that if
$\omega \in \Lambda^p(M)$ then
\begin{equation} \label{tau}
\mu(\omega) = i^{\frac{n(n-1)}{2}} \: (-1)^{\frac{p(p-1)}{2}} * (\omega).
\end{equation}  
One can check that $\mu^2 = 1$ and $\mu (d + d^*) + (d + d^*) \mu = 0$.

Clearly the operator $(d + d^*) \otimes \Id_N$, 
acting on $\Lambda^*(M) \otimes \C^N$, is also invertible.
Consider the $u(N)$-valued $1$-form 
$F^* \omega_\epsilon$ on $M$. As $F$ is distance-nonincreasing,
\begin{equation}
\parallel F^* \omega_\epsilon \parallel_\infty \: \le \: 
\parallel \omega_\epsilon \parallel_\infty \: = \epsilon
\parallel \omega \parallel_\infty.
\end{equation}  
Let $e(F^* \omega_\epsilon)$ denote exterior multiplication by 
$F^* \omega_\epsilon$, acting on $\Lambda^*(M) \otimes \C^N$ and let
$i(F^* \omega_\epsilon)$ denote interior multiplication by 
$F^* \omega_\epsilon$. By making $\epsilon$ small enough, the operator
$e(F^* \omega_\epsilon) - i(F^* \omega_\epsilon)$ has arbitrarily small
norm and so the operator 
$\left( (d + d^*) \otimes \Id_N \right) + 
e(F^* \omega_\epsilon) - i(F^* \omega_\epsilon)$ is also invertible.

Put $D = \left( d \otimes \Id_N \right) + 
e(F^* \omega_\epsilon)$. Then $D$ is exterior differentiation, using the
connection $F^* \omega_\epsilon$, and 
\begin{equation}
D + D^* = \left( (d + d^*) \otimes \Id_N \right) + 
e(F^* \omega_\epsilon) - i(F^* \omega_\epsilon).
\end{equation}
As $(d + d^*) \otimes \Id_N$ and $D + D^*$ anticommute with
$\mu \otimes \Id_N$, they have well-defined indices which happen to vanish,
as the operators are invertible. On the other hand, let $L(M)$
be the Hirzebruch $L$-form. The relative index
theorem of Gromov and Lawson \cite{Donnelly (1987),Gromov-Lawson (1983)} 
says that
\begin{equation}
\ind(D + D^*) - \ind( (d + d^*) \otimes \Id_N ) =
\int_M L(M) \wedge 
\left[ \Tr \left( e^{- \frac{F^* \Omega_\epsilon}{2\pi i}} \right) 
- N \right].
\end{equation}
As $F$ is proper, the de Rham cohomology class of $\Tr \left( 
e^{- \frac{F^* \Omega_\epsilon}{2\pi i}} \right) - N  = F^* \left[ \Tr \left( 
e^{- \frac{\Omega_\epsilon}{2\pi i}} \right) - N \right]$  
is well-defined as
a multiple of the fundamental class $[M] \in \HH^n_c(M; \R)$.
As the series for $L(M)$ starts off as 
$L(M) = 1 + \ldots$, we obtain
\begin{align}
\ind(D + D^*) - \ind( (d + d^*) \otimes \Id_N ) & =
\int_M \left[ \Tr \left( 
e^{- \frac{F^* \Omega_\epsilon}{2\pi i}} \right) - N \right] \\
& =
\int_M F^* \left[
\Tr \left( e^{- \frac{\Omega_\epsilon}{2\pi i}} \right) - N \right] \notag \\
& = \deg(F) \int_{\R^n} \left[
\Tr \left( e^{- \frac{\Omega_\epsilon}{2\pi i}} \right) - N \right]
\neq 0 \notag.
\end{align}
This contradicts the vanishing of 
$\ind(D + D^*)$ and $\ind( (d + d^*) \otimes \Id_N )$. Thus zero must be
in the spectrum of $M$ after all.

Now suppose that $n$ is odd. As $M$ is hyperEuclidean, so is $\R \times M$.
With respect to the decomposition $\Lambda^*(\R \times M) = \Lambda^*(\R)
\otimes \Lambda^*(M)$, the Laplace-Beltrami operator on $\R \times M$
decomposes as
\begin{equation}
\triangle_{\R \times M} = \left( \triangle_{\R} \otimes I \right) +
\left( I \otimes \triangle_{M} \right).
\end{equation}
Then
\begin{equation}
\sigma(\triangle_{\R \times M}) = \{ \lambda_1 + \lambda_2 : 
\lambda_1 \in [0, \infty) \text{ and } \lambda_2 \in \sigma(\triangle_M) \}.
\end{equation}
From what has already been proved, $0 \in \sigma(\triangle_{\R \times M})$.
It follows that $0 \in \sigma(\triangle_M)$. 
\end{pf} 
\noindent
{\bf Remarks : }\\
1. We have shown that if $M$ is hyperEuclidean then $0 \in \sigma(\triangle_p)$
for some $p$. One can ask whether the number $p$ can be pinned down.  In 
general, when computing the index of the operator $d + d^*$, the 
differential forms outside of the middle dimensions do not contribute.
This is a reflection of the fact that the signature of a closed manifold
can be computed using only the middle-dimensional cohomology. So this
gives some reason to think that if $\dim(M)$ is even then
$0 \in \sigma \left(\triangle_{\frac{dim(M)}{2}}\right)$.

Unfortunately, the operator $(D + D^*)^2$ does not preserve the degree of
a differential form and so we cannot use the above proof to 
reach the desired conclusion.  However, with a more refined
index theorem \cite[Theorem 6.24]{Roe (1993)}, one can indeed conclude that
$0 \in \sigma \left(\triangle_{\frac{dim(M)}{2}}\right)$ 
if $\dim(M)$ is even and that
$0 \in \sigma \left(\triangle_{\frac{dim(M) \pm 1}{2}} \right)$ 
if $\dim(M)$ is odd.\\
2. If $M$ is an irreducible noncompact globally symmetric space
$G/K$, with $G = \Isom(M)$ and $K$ a maximal compact subgroup,  then
one can say more about the bottom of the spectrum. If
$\rk(G) = \rk(K)$ then $\Ker \left(\triangle_{\frac{dim(M)}{2}} \right)$ 
is infinite-dimensional
and the spectrum of $\triangle$ is bounded away from zero otherwise.
If $\rk(G) > \rk(K)$ then $\Ker(\triangle) = 0$ and $0 \in \sigma(\triangle_p)$
if and only if $p \in \left[ \frac{dim(M)}{2} - \frac{rk(G) - rk(K)}{2},
\frac{dim(M)}{2} + \frac{rk(G) - rk(K)}{2} \right]$ 
\cite[Section VIIB]{Lott (1990)}. \\

Finally, we state a result about uniformly contractible Riemannian manifolds.
\begin{definition} \cite[p. 29]{Gromov (1993)}
A metric space $Z$ has finite asymptotic dimension if there
is an integer $n$ such that for any $r > 0$, there is a covering 
$Z = \bigcup_{i \in I} C_i$ 
of $Z$ by subsets of uniformly bounded diameter so
that no metric ball of radius $r$ in $Z$ intersects more than $n+1$
elements of $\{C_i\}_{i \in I}$. The smallest such integer $n$ is called
the asymptotic dimension $\asdim_+(Z)$ of $Z$.
\end{definition} 

\begin{proposition} (Yu \cite{Yu (1995)}) \label{Yu}
If $M$ is a uniformly contractible
Riemannian manifold with finite asymptotic dimension then $0 \in 
\sigma(\triangle_p)$ for some $p$.
\end{proposition}

The proof of Proposition \ref{Yu} uses methods of coarse index theory
\cite{Roe (1993)}. 
\section{Very Low Dimensions}

In this section we show that the answer to the zero-in-the-spectrum 
question is ``yes'' for one-dimensional simplicial complexes and 
two-dimensional Riemannian manifolds.
\subsection{One Dimension}

As a one-dimensional manifold is either $S^1$ or $\R$, zero is clearly in the
spectrum.

A more interesting problem is to consider a connected
one-dimensional simplicial
complex $K$. Let $V$ be the set of vertices of $K$ and let $E$ be the
set of oriented edges of $K$. That is, an element $e$ of $E$ consists of
an edge of $K$ and an ordering $(s_e, t_e)$ of $\partial e$. We
let $- e$ denote the same edge with the reverse ordering of $\partial e$. 
For $x \in V$,
let $m_x$ denote the number of unoriented edges of which $x$ is a boundary.
We assume that $m_x < \infty$ for all $x$. Put
\begin{align}
C^0(K) &= \{ f : V \rightarrow \C \text{ such that }
\sum_{x \in V} \: m_x \: |f(x)|^2 < \infty\},\\
C^1(K) & = \{ F : E \rightarrow \C \text{ such that }
F(-e)=-F(e) \text{ and } \frac12 \sum_{e \in E} |F(e)|^2 < \infty\}. \notag
\end{align}
Then $C^0(K)$ and $C^1(K)$ are Hilbert spaces. The edges $e$ such that
$s_e = t_e$ do not enter in $C^1(K)$ and can be deleted for our purposes.
The weighting used to define $C^0(K)$ is natural in certain
respects \cite{Dodziuk-Kendall (1986)}. 

There is a bounded operator
$d : C^0(K) \rightarrow C^1(K)$ given by
$(df)(e) = f(t_e) - f(s_e)$. Define the Laplace-Beltrami operators by 
$\triangle_0 = d^* d$ and $\triangle_1 = d d^*$.
An element of $\Ker(\triangle_1)$ is an $F \in C^1(K)$ such that
for each vertex $x$ the total current flowing into $x$ vanishes, i.e.
$\sum_{e \in E : t_e = x} F(e) = 0$.

The next proposition is essentially due to Gromov \cite[p. 236]{Gromov (1993)},
who proved it in the case when $\{m_x\}_{x \in V}$ is bounded.
\begin{proposition} \label{simplicial}
$0 \in \sigma(\triangle_0)$ or $0 \in \sigma(\triangle_1)$.
\end{proposition}
\begin{pf}
As the nonzero spectra of $d^*d$ and $d d^*$ are the same, for our purposes
it suffices to consider $\sigma(\triangle_0)$ and $\Ker(\triangle_1)$.
We argue by contradiction.
Suppose that $0 \notin \sigma(\triangle_0)$ and $\Ker(\triangle_1) = 0$.
First, $K$ must be infinite, as otherwise $\Ker(\triangle_0) \ne 0$.
Second, $K$ must be a tree, as if $K$ had a loop then we could create a nonzero
element of $\Ker(\triangle_1)$ by letting a current of unit strength
flow around the loop.

We now show that $K$ has lots of branching.
For $x, y \in V$, let $[x,y]$ be the geodesic arc from $x$ to $y$ and let
$(x,y)$ be its interior. Let $d(x,y)$ be the number of edges in $[x,y]$. 
\begin{lemma} There is a constant $L > 0$ such that if $d(x,y) > L$ then
there is an infinite subtree of $K$ which intersects $(x,y)$ but does
not contain $x$ or $y$.
\end{lemma}
\begin{pf}
If the lemma is not true
then for all $N > 1$, there are vertices $x$ and $y$ such
that $d(x,y) > N$ but  
there are no infinite subtrees as in the statement of the lemma.  In
other words, the connected component $C$ of $K - \{x\} - \{y\}$ which
contains $(x,y)$ is finite. As $K$ is a tree, $x$ is only connected to
the vertices in $C$ by a single edge, and similarly for $y$ (see Fig. 5).
Define $f \in C^0(K)$ by
\begin{equation}
f(v) = 
\begin{cases}
1 & \text{ if } v \in C,\\
0 & \text{ otherwise.}
\end{cases}
\end{equation}
Then
\begin{equation}
\frac{\langle df, df \rangle}{\langle f,f \rangle} \le 
\frac{2}{2(d(x,y) - 1)} \le 
\frac{1}{N}.
\end{equation}
As $N$ can be taken arbitrarily large, 
this contradicts the assumption that $0 \notin \sigma(\triangle_0)$.
\end{pf} 

It follows that $K$ contains a subtree $K^\prime$ which is topologically
equivalent to an infinite triadic tree, with the distances between
branchings at most $L$ (see Fig. 6). We can create a nonzero 
square-integrable harmonic
$1$-cochain $F^\prime$ on $K^\prime$ by letting a unit current flow 
through it, as in Fig. 6. Let $F \in C^1(K)$ be the extension of $F^\prime$
by zero to $K$. If $x$ is a vertex of $K^\prime$ then the total current
flowing into $x$ is still zero, as no new current comes in along the edges of
$K - K^\prime$.   Thus $\Ker(\triangle_1) \neq 0$, which is a contradiction. 
\end{pf}

\subsection{Two Dimensions}
\begin{proposition} (Lott, Dodziuk) \label{surfaces}
The answer to the zero-in-the-spectrum question is ``yes'' if $M$ is
a two-dimensional manifold.
\end{proposition}
\begin{pf}
The Hodge decomposition gives
\begin{align}
\Lambda^0(M) & = \Ker(\triangle_0) \oplus \Lambda^0(M)/\Ker(d), \\
\Lambda^1(M) & = \Ker(\triangle_1) \oplus \overline{d \Lambda^0(M)} 
\oplus * \overline{d \Lambda^0(M)},\\
\Lambda^2(M) & = * \Ker(\triangle_0) \oplus * (\Lambda^0(M)/\Ker(d)).
\end{align}
Thus it is enough to look at $\Ker(\triangle_0)$, $\Ker(\triangle_1)$ and
$\sigma \left(\triangle \text{ on } \Lambda^0(M)/\Ker(d) \right)$.

We argue by contradiction.  Assume that zero is not in the spectrum. 
By Proposition \ref{compactsup}, $\Image(\HH^1_c(M) \rightarrow \HH^1(M)) = 0$.
Thus $M$ must be planar, in the sense of either of
the following two equivalent conditions
: \\
1. Any simple closed curve in $M$ separates it into two pieces.\\
2. $M$ is diffeomorphic to the complement of a closed subset of $S^2$.\\

As $\Ker(\triangle_0) = 0$, $M$ cannot be $S^2$. 
By Proposition \ref{conformal}, the possible existence of nonzero
square-integrable harmonic $1$-forms on $M$ only depends on the
underlying Riemann surface coming from the Riemannian metric
on $M$. 

We recall some notions from Riemann surface theory
\cite{Ahlfors-Sario (1960)}. A function
$f \in C^\infty(M)$ is {\em superharmonic} if $\triangle_0 f > 0$. 
(This is a conformally-invariant statement.)
The Riemann surface underlying $M$ is {\em hyperbolic} 
if it has a positive superharmonic function and
{\em parabolic} otherwise. If $M$ is planar and hyperbolic then
there is a nonconstant harmonic function $f \in C^\infty(M)$ such that
$\int_M df \wedge *df < \infty$ \cite[p. 208]{Ahlfors-Sario (1960)}. 
Then $df$ would be a nonzero element of $\Ker(\triangle_1)$. Thus
$M$ must be parabolic.

Put $\lambda_0 = \inf(\sigma(\triangle_0))$. Choose some $\lambda$ such that
$0 < \lambda < \lambda_0$. Then there is a positive 
$f \in C^\infty(M)$ (not square-integrable!)
such that $\triangle_0 f = \lambda f$
\cite[Theorem 2.1]{Sullivan (1987)}. However, this contradicts the
parabolicity of $M$.
\end{pf}

We do not know of any result analogous to Proposition \ref{surfaces} for
general two-dimensional simplicial complexes, say uniformly finite. See, 
however, Subsection \ref{Two and Three Dimensions}.
\section{Universal Covers}

Suppose that $M$ is the universal cover of
a compact Riemannian manifold $X$. We give $M$ the pulled-back Riemannian
metric and consider the Laplace-Beltrami operator
$\triangle_p$ on $M$. There are numerical invariants which measure
the density of $\sigma(\triangle_p)$ near zero, the so-called $L^2$-Betti
numbers $\{b^{(2)}_p(X)\}_{p \ge 0}$ 
and Novikov-Shubin invariants $\{\alpha_{p+1}(X)\}_{p \ge 0}$.
We refer to \cite{Lott-Lueck (1995),Lueck (1996),Pansu (1996)} 
for the definitions of these invariants.  We will
only need the following properties :\\ \\
{\bf Properties :}
1. $b^{(2)}_p(X) = 0 \iff \Ker(\triangle_p) = 0$.\\
2. $0 \notin \sigma(\triangle_p \text{ on } \Lambda^p(M)/\Ker(d)) \iff
\alpha_{p+1} = \infty^+$.\\
3. $b^{(2)}_p(X)$ and $\alpha_p(X)$ are homotopy-invariants of $X$.\\
4. $b^{(2)}_0(X)$, $b^{(2)}_1(X)$, 
$\alpha_1(X)$ and $\alpha_2(X)$ only depend on
$\pi_1(X)$.\\
5. $b^{(2)}_0(X) = 0$ if and only if $\pi_1(X)$ is infinite.\\
6. $\alpha_1(X) = \infty^+$ if and only if $\pi_1(X)$ is 
finite or nonamenable.\\
7. The Euler characteristic of $X$ satisfies
\begin{equation}
\chi(X) = \sum_p (-1)^p \: b^{(2)}_p(X)
\end{equation}
\noindent
8. If $X^n$ is closed then $b^{(2)}_{n-p}(X) = b^{(2)}_p(X)$.\\
9. If $X^{4k}$ is closed 
then there are nonnegative numbers
$b^{(2)}_{2k,\pm}(X)$ such that $b^{(2)}_{2k}(X) = b^{(2)}_{2k,+}(X)
+ b^{(2)}_{2k,-}(X)$ and the signature of $X$ satisfies
\begin{equation}
\tau(X) = b^{(2)}_{2k,+}(X) - b^{(2)}_{2k,-}(X).
\end{equation}

One can extend properties 1-7 from compact Riemannian manifolds $X$ to
finite $CW$-complexes $K$.

In what follows, $\Gamma$ will denote a finitely-presented group.
Given a presentation of $\Gamma$, there is an associated $2$-dimensional
$CW$-complex $K$ which we call the {\em presentation complex}. To form it,
make a bouquet of
circles indexed by the generators of $\Gamma$. Attach $2$-cells based on
the relations of $\Gamma$. (We allow trivial or repeated relations
in the presentation.) This is the presentation complex.

\begin{definition} \label{group}
Put $b_0^{(2)}(\Gamma) = b_0^{(2)}(K)$,
$b_1^{(2)}(\Gamma) = b_1^{(2)}(K)$, $\alpha_1(\Gamma) = \alpha_1(K)$
and $\alpha_2(\Gamma) = \alpha_2(K)$.
\end{definition}

By Property 4 above, Definition \ref{group} makes sense in that
the choice of presentation of $\Gamma$ does not matter.

\subsection{Big and Small Groups}
\begin{definition}
The group $\Gamma$ is {\em big} if it is nonamenable, $b_1^{(2)}(\Gamma) = 0$
and $\alpha_2(\Gamma) = \infty^+$. Otherwise, $\Gamma$ is {\em small}.
\end{definition}
\begin{proposition} \label{smallspec}
Let $X$ and $M$ be as above.
The group $\pi_1(X)$ is small if and only if $0 \in \sigma(\triangle_0)$ or
$0 \in \sigma(\triangle_1)$.
\end{proposition}
\begin{pf}
This follows immediately from Properties 1, 2, 4, 5 and 6 above.
\end{pf}

The question arises as to which groups are big and which are small.
Clearly any amenable group is small.
\begin{proposition} \label{surfacegp}
Fundamental groups of compact surfaces are small.
\end{proposition}
\begin{pf}
Suppose that $\Sigma$  is a compact surface and $\Gamma =
\pi_1(\Sigma)$. 
If $\Sigma$ has boundary then $\Gamma$ is a free group $F_j$ on some
number $j$ of generators. If $j = 0$ or $j =1$ then $\Gamma$ is amenable.
If $j > 1$ then $b_1^{(2)}(\Gamma) = j - 1 > 0$.

Suppose now that $\Sigma$ is closed. If $\chi(\Sigma) \ge 0$ then
$\Gamma$ is amenable. If $\chi(\Sigma) < 0$ then
$b_1^{(2)}(\Gamma) = - \chi(\Sigma) > 0$.  
\end{pf}

We now extend Proposition \ref{surfacegp} to $3$-manifold groups. We
use some facts about compact connected $3$-manifolds
$Y$, possibly with boundary. (See, for example, 
\cite[Section 6]{Lott-Lueck (1995)}). Again, all of our manifolds are
assumed to be oriented. 
First, $Y$ has a decomposition as a connected sum $Y = Y_1 \# Y_2 \# \ldots
\# Y_r$ of {\em prime} $3$-manifolds. A prime $3$-manifold
is {\em exceptional} if it is closed and no finite cover of it is 
homotopy-equivalent to a Seifert, Haken or hyperbolic $3$-manifold.
No exceptional prime $3$-manifolds are known and it is likely that
there are none.
\begin{proposition} (Lott-L\"uck) \label{3mfldgp}
Suppose that $Y$ is a compact connected oriented $3$-manifold, possibly
with boundary, none of whose prime factors are exceptional.  Then
$\pi_1(Y)$ is small.
\end{proposition}
\begin{pf}
We argue by contradiction. Suppose that $\pi_1(Y)$ is big.
First, $\pi_1(Y)$ must be infinite.
If $\partial Y$ has any connected components which are
$2$-spheres then we can cap them
off with $3$-balls without changing $\pi_1(Y)$. So we can assume that
$\partial Y$ does not have any $2$-sphere components. In particular,
$\chi(Y) = \frac12 \chi(\partial Y) \le 0$.
From \cite[Theorem 0.1.1]{Lott-Lueck (1995)},
\begin{equation}
b_1^{(2)}(Y) = (r-1) - \sum_{i=1}^r \frac{1}{|\pi_1(Y_i)|} - \chi(Y).
\end{equation}
As this must vanish, we have $\chi(Y) = 0$ and either\\
1. $\{ | \pi_1(Y_i)| \}_{i=1}^r = \{2,2,1,\ldots,1\}$ or\\
2. $\{ | \pi_1(Y_i)| \}_{i=1}^r = \{\infty,1,\ldots,1\}$.

It follows that $\partial Y$ is empty or a disjoint union of $2$-tori.
As there are no $2$-spheres in $\partial Y$, if $|\pi_1(Y_i)| = 1$
then $Y_i$ is a homotopy $3$-sphere. Thus $Y$ is homotopy-equivalent
to either\\
1. $\R P^3 \# \R P^3$ or\\
2. A prime $3$-manifold $Y^\prime$ 
with infinite fundamental group whose boundary
is empty or a disjoint union of $2$-tori.

If $Y$ is homotopy-equivalent to $\R P^3 \# \R P^3$ then $\pi_1(Y)$ is
amenable, which is a contradiction.  So we must be
in the second case.  Using Property 3, we may assume that $Y = Y^\prime$.
Then as $Y$ is prime, it follows from 
\cite[Chapter 1]{McCullough-Miller (1986)} that
either $Y = S^1 \times D^2$ or $Y$ has incompressible (or empty) boundary.
If $Y = S^1 \times D^2$ then $\pi_1(Y)$ is amenable. If
$Y$ has incompressible (or empty) boundary then
from \cite[Theorem 0.1.5]{Lott-Lueck (1995)}, $\alpha_2(Y) \le 2$ unless
$Y$ is a closed $3$-manifold with an $\R^3$, $\R \times S^2$ or $Sol$
geometric structure.  In the latter cases, $\Gamma$ is amenable. 
Thus in any case, we get a contradiction. 
\end{pf}

The next proposition gives examples of big groups.
\begin{proposition}
1. A product of two nonamenable groups is big.\\
2. If $Y$ is a closed nonpositively-curved locally symmetric space of
dimension greater than three, with no Euclidean factors in $\widetilde{Y}$,
then $\pi_1(Y)$ is big.
\end{proposition}
\begin{pf}
1. Suppose that 
$\Gamma = \Gamma_1 \times \Gamma_2$ with $\Gamma_1$ and $\Gamma_2$
nonamenable. Then $\Gamma$ is nonamenable. 
Let $K_1$ and $K_2$ be presentation complexes with fundamental
groups $\Gamma_1$ and $\Gamma_2$, respectively. Put $K = K_1 \times K_2$.
Then $\Gamma = 
\pi_1(K)$. Let $\triangle_p(\widetilde{K})$, 
$\triangle_p(\widetilde{K_1})$ and
$\triangle_p(\widetilde{K_2})$ 
denote the Laplace-Beltrami operator on $p$-cochains on
$\widetilde{K}$, $\widetilde{K_1}$ and $\widetilde{K_2}$, respectively,
as defined in Subsection \ref{Two and Three Dimensions} below. Then
\begin{align} \label{Kunneth}
\inf(\sigma(\triangle_1(\widetilde{K}))) = & \min \left(
 \inf(\sigma(\triangle_1(\widetilde{K_1}))) + 
\inf(\sigma(\triangle_0(\widetilde{K_2}))), \right. \\
& \left. \hspace{.38in} \inf(\sigma(\triangle_0(\widetilde{K_1}))) +
\inf(\sigma(\triangle_1(\widetilde{K_2}))) \right) > 0. \notag
\end{align}
Using Proposition \ref{smallspec}, the first part of the proposition follows.\\
2. If $\widetilde{Y}$ is irreducible then part 2. of the
proposition follows from the second remark after Proposition \ref{hyper}.
If $\widetilde{Y}$ is reducible then we can use an argument similar to
(\ref{Kunneth}). \end{pf}
\subsection{Two and Three Dimensions} \label{Two and Three Dimensions}

In this subsection we relate the zero-in-the-spectrum question to a
question in combinatorial group theory.
Let $K$ be a finite connected $2$-dimensional $CW$-complex.
Let $\widetilde{K}$ be its universal cover. Let $C^*(\widetilde{K})$ denote the
Hilbert space of square-integrable
cellular cochains on $\widetilde{K}$. There is a cochain complex
\begin{equation}
0 \longrightarrow C^0(\widetilde{K}) \stackrel{d_0}{\longrightarrow}
C^1(\widetilde{K}) \stackrel{d_1}{\longrightarrow}
C^2(\widetilde{K}) \longrightarrow 0.
\end{equation}
Define the Laplace-Beltrami operators by
$\triangle_0 = d_0^* d_0$, $\triangle_1 = d_0 d_0^* + d_1^* d_1$ and
$\triangle_2 = d_1 d_1^*$. These are bounded self-adjoint operators and 
so we can talk about zero being in the
spectrum of $\widetilde{K}$.
\begin{proposition} \label{groups}
Zero is not in the spectrum of $\widetilde{K}$ if and only if
$\pi_1(K)$ is big and $\chi(K) = 0$.
\end{proposition}
\begin{pf}
Suppose that zero is not in the spectrum of $\widetilde{K}$. From the
analog of Proposition \ref{smallspec}, $\Gamma$ must be big.
Furthermore, from Properties 1 and 7, $\chi(K) = 0$.

Now suppose that $\pi_1(K)$ is big and $\chi(K) = 0$.
From the analog of Proposition \ref{smallspec}, $0 \notin \sigma(\triangle_0)$
and $0 \notin \sigma(\triangle_1)$. In particular, $\Ker(\triangle_0) =
\Ker(\triangle_1) = 0$. From Properties 1 and 7, $\Ker(\triangle_2) = 0$.
As $C^2(\widetilde{K}) = \Ker(\triangle_2) \oplus 
\overline{d_1 C^1(\widetilde{K})}$, we conclude that 
$0 \notin \sigma(\triangle_2)$.
\end{pf}

Let $\Gamma$ be a finitely-presented group. Consider a fixed presentation
of $\Gamma$ consisting of $g$ generators and $r$ relations. Let $K$ be
the corresponding presentation complex. Then $\chi(K) = 1 - g + r$.
Thus zero is not in the spectrum of $\widetilde{K}$ if and only if
$\pi_1(K)$ is big and $g - r = 1$.

Recall that the {\em deficiency} $\Def(\Gamma)$
is defined to be the
maximum, over all finite presentations of $\Gamma$, of $g - r$. 
If $b_1^{(2)}(\Gamma) = 0$ then from the equation
\begin{equation}
\chi(K) = 1 - g + r = b_0^{(2)}(\Gamma) - b_1^{(2)}(\Gamma) + 
b_2^{(2)}(K),
\end{equation}
we obtain $\Def(\Gamma) \le 1$.  
This is the case, for example, when $\Gamma$ is big or
when $\Gamma$ is amenable \cite{Cheeger-Gromov (1986)}.

As any finite connected $2$-dimensional $CW$-complex is homotopy-equivalent to
a presentation complex, it follows from Proposition \ref{groups} that
the answer to the zero-in-the-spectrum question is ``yes'' for
universal covers of such complexes if and only if the following
conjecture is true.

\begin{conjecture}
If $\Gamma$ is a big group then $\Def(\Gamma) \le 0$.
\end{conjecture}
\noindent
{\bf Remark : } If $\pi_1(K)$ has property $T$ then the ordinary
first Betti number of $K$ vanishes, and so
$\chi(K) = 1 + b_2(K) > 0$. Thus zero lies in the spectrum of 
$\widetilde{K}$. \\

Now let $Y$ be a $3$-manifold satisfying the conditions of Proposition
\ref{3mfldgp}. If $\partial Y \neq \emptyset$, we define $\triangle_p$ on
$\widetilde{Y}$ using absolute boundary conditions on 
$\partial \widetilde{Y}$.
\begin{proposition}
Zero lies in the spectrum of $\widetilde{Y}$.
\end{proposition}
\begin{pf}
This is a corollary of Propositions \ref{smallspec} and \ref{3mfldgp}.
\end{pf} 

\subsection{Four Dimensions}
In the subsection we relate the zero-in-the-spectrum question to a 
question about Euler characteristics of closed $4$-dimensional manifolds.

If $M$ is a Riemannian $4$-manifold then the Hodge decomposition gives
\begin{align} \label{4hodge}
\Lambda^0(M) & = \Ker(\triangle_0) \oplus \Lambda^0(M)/\Ker(d), \\
\Lambda^1(M) & = \Ker(\triangle_1) \oplus \overline{d \Lambda^0(M)} 
\oplus \Lambda^1(M)/\Ker(d), \notag \\
\Lambda^2(M) & = \Ker(\triangle_2) \oplus \overline{d\Lambda^1(M)}
\oplus * \overline{d\Lambda^1(M)}, \notag \\
\Lambda^3(M) & = * \Ker(\triangle_1) \oplus * \overline{d \Lambda^0(M)}
\oplus *(\Lambda^1(M)/\Ker(d)), \notag  \\
\Lambda^4(M) & = * \Ker(\triangle_0) \oplus * (\Lambda^0(M)/\Ker(d)).
\notag 
\end{align}
Thus for the zero-in-the-spectrum question, it is enough to consider
$\Ker(\triangle_0)$, $\Ker(\triangle_1)$,
$\sigma(\triangle_0 \text{ on } \Lambda^0/\Ker(d))$,
$\sigma(\triangle_1 \text{ on } \Lambda^1/\Ker(d))$ and
$\Ker(\triangle_2)$.

Let $\Gamma$ be a finitely-presented group.
Recall that $\Gamma$ is the fundamental
group of some 
closed $4$-manifold. To see this, take a finite presentation of
$\Gamma$.  Embed the resulting presentation complex
in $\R^5$ and take the boundary of a regular neighborhood to get the manifold.

Now consider the Euler
characteristics of all 
closed $4$-manifolds $X$ with fundamental group $\Gamma$.
Given $X$, we have $\chi(X \# \C P^2) = \chi(X) + 1$. Thus it is easy
to make the Euler characteristic big.  However, it is not so easy to make
it small. From what has been said,
\begin{equation}
\{ \chi(X) : X \text{ is a closed 
connected oriented 4-manifold with }
\pi_1(X) = \Gamma \} = \{n \in \Z : n \ge q(\Gamma)\}
\end{equation}
for some $q(\Gamma)$. {\em A priori}
$q(\Gamma) \in \Z \cup \{ - \infty\}$, but in fact $q(\Gamma) \in \Z$
\cite[Th\'eor\`eme 1]{Hausmann-Weinberger (1985)}.
(This also follows from (\ref{Euler}) below.) It is a basic problem
in $4$-manifold topology to get good estimates of $q(\Gamma)$.

Suppose that $\pi_1(X) = \Gamma$. From Properties 4, 7 and 8 above,
\begin{equation} \label{Euler}
\chi(X) = 2 \: b^{(2)}_0(\Gamma) - 2 \: b^{(2)}_1(\Gamma) + b^{(2)}_2(X).
\end{equation}
In particular, if $b^{(2)}_1(\Gamma) = 0$ then $\chi(X) \ge 0$ and so
$q(\Gamma) \ge 0$. This is the case, for example, when $\Gamma$ is big or
when $\Gamma$ is amenable \cite{Cheeger-Gromov (1986)}.

\begin{proposition}
Let $X$ be a closed $4$-manifold. Then zero is not in the spectrum of
$\widetilde{X}$ if and only if $\pi_1(X)$ is big and $\chi(X) = 0$.
\end{proposition}
\begin{pf}
Suppose that zero is not in the spectrum of $\widetilde{X}$.
Then from Proposition
\ref{smallspec}, $\pi_1(X)$ must be big.  Furthermore, 
$\Ker(\triangle_2) = 0$. From Property 1 and (\ref{Euler}), $\chi(X) = 0$.

Now suppose that $\pi_1(X)$ is big and $\chi(X) = 0$. From
Proposition \ref{smallspec},  
$0 \notin \sigma(\triangle_0)$ and $0 \notin \sigma(\triangle_1)$. 
From Property 1 and 
(\ref{Euler}), $\Ker(\triangle_2) = 0$. Then from
(\ref{4hodge}), zero is not in the spectrum of
$\widetilde{X}$.
\end{pf}
{\bf Remark : } If zero is not in the spectrum of $\widetilde{X}$ then it
follows from Property 9 that in addition, $\tau(X) = 0$. Also, if $\pi_1(X)$
satisfies the Strong Novikov Conjecture then it will follow from
Corollary \ref{SNC} that $\nu_*([X])$ vanishes in $\HH_4(B\pi_1(X); \C)$.\\

In summary, we have shown that the answer to the zero-in-the-spectrum question
is ``yes'' for universal covers of closed $4$-manifolds if and only if the
following conjecture is true.
\begin{conjecture} \label{4conj}
If $\Gamma$ is a big group then $q(\Gamma) > 0$.
\end{conjecture}

We now give some partial positive results on the 
zero-in-the-spectrum question for universal covers of closed $4$-manifolds.
Recall that there is a notion, due to Thurston, 
of a manifold having a geometric structure. This is especially important
for $3$-manifolds. The $4$-manifolds with geometric structures have been
studied by Wall \cite{Wall (1986)}.

\begin{proposition}
Let $X$ be a closed $4$-manifold. Then zero is in the spectrum of 
$\widetilde{X}$ if\\
1. $\pi_1(X)$ has Property $T$ or\\
2. $X$ has a geometric structure (and an arbitrary Riemannian metric) or\\
3. $X$ has a complex structure (and an arbitrary Riemannian metric).
\end{proposition}
\begin{pf}
1. If $X$ has Property $T$ then the ordinary first Betti number of $X$ 
vanishes. Thus $\chi(X) = 2 + b_2(X) > 0$. Part 1. of the proposition follows.
\\
2. The geometries of \cite{Wall (1986)} all fall into at least one of
the following classes :\\
a. $b_0^{(2)} \neq 0$ : $S^4$, $S^2 \times S^2$, $\C P^2$.\\
b. $0 \in \sigma(\triangle_0 \text{ on } \Lambda^0/\Ker(d))$ :
$\R^4$, $S^3 \times \R$, $S^2 \times \R^2$, $Nil^3 \times \R$, 
$Nil^4$, $Sol^4_0$, $Sol^4_1$, $Sol^4_{m,n}$.\\
c. $b_1^{(2)} \neq 0$ : $S^2 \times H^2$.\\
d. $0 \in \sigma(\triangle_1 \text{ on } \Lambda^1/\Ker(d))$ :
$H^3 \times \R$, $\widetilde{SL_2} \times \R$, $H^2 \times \R^2$.\\
e. $\chi > 0$ :  $H^4$, $H^2 \times H^2$, $\C H^2$. \\
Part 2. of the proposition follows.\\
3. Suppose that zero is not in the spectrum of $\widetilde{X}$.
From Properties 7 and 9, 
$\chi(X) = \tau(X) = 0$. From the classification of complex surfaces,
$X$ has a geometric structure \cite[p. 148-149]{Wall (1986)}. This contradicts
part 2. of the proposition.
\end{pf}

\subsection{More Dimensions}
In this subsection we give some partial positive results about the
zero-in-the-spectrum question for covers of compact manifolds of
arbitrary dimension.
Let us first recall some facts about
index theory \cite{Kasparov (1988)}.  
Let $X$ be a closed Riemannian manifold.
If $\dim(X)$ is even, consider the
operator $d + d^*$ on $\Lambda^*(X)$. Give $\Lambda^*(X)$ the $\Z_2$-grading
coming from (\ref{tau}). Then the signature $\tau(X)$ equals the
index of $d+d^*$. To say this in a more complicated way,
the operator $d + d^*$ defines a element $[d + d^*]$ of the K-homology group
$K_0(X)$. Let $\nu : X \rightarrow \point$ be the (only) map from $X$
to a point. Then $\nu_*([d+d^*]) \in K_0(\point)$. There is a map
$A : K_0(\point) \rightarrow K_0(\C)$ which is the identity, as both
sides are $\Z$. So we can say that $\tau(X) = A(\nu([d+d^*])) \in K_0(\C)$.

Now let $M$ be a normal cover of $X$ with covering group $\Gamma$.
The fiber bundle $M \rightarrow X$ is classified by a map
$\nu : X \rightarrow B\Gamma$, defined up to homotopy. Let $\widetilde{d}$ be
exterior differentiation on $M$.
Consider the operator $\widetilde{d} + \widetilde{d}^*$. 
Taking into account the action of
$\Gamma$ on $M$, one can define a refined index 
$\ind(\widetilde{d} + \widetilde{d}^*) \in K_0(C^*_r \Gamma)$, where
$C^*_r \Gamma$ is the reduced group $C^*$-algebra of $\Gamma$.
  
We recall the statement of the Strong Novikov
Conjecture (SNC). 
This is a conjecture about a countable discrete group $\Gamma$, namely that
the assembly map $A : K_*(B\Gamma) \rightarrow K_*(C^*_r \Gamma)$
is rationally injective. Many groups of a geometric origin, such as
discrete subgroups of connected
Lie groups or Gromov-hyperbolic groups, are known to
satisfy SNC.  There are no known groups which do not satisfy SNC.
\begin{proposition} \label{SNC1}
Let $X$ be a closed Riemannian
manifold with a surjective homomorphism $\pi_1(X) \rightarrow \Gamma$.
Let $M$ be the induced normal $\Gamma$-cover of $X$.
Suppose that $\Gamma$ satisfies SNC.
Let $L(X) \in \HH^{*}(X; \C)$ be the Hirzebruch $L$-class of $X$ and let
$*L(X) \in \HH_*(X; \C)$ be its Poincar\'e dual.
Then if
$\nu_*(*L(X)) \neq 0$ in $H_*(B\Gamma; \C)$, zero lies in the spectrum of
$M$. In fact,
$0 \in \sigma \left(\triangle_{\frac{dim(X)}{2}}\right)$ 
if $\dim(X)$ is even and 
$0 \in \sigma \left(\triangle_{\frac{dim(X) \pm 1}{2}} \right)$ 
if $\dim(X)$ is odd.
\end{proposition}
\begin{pf}
Suppose first that $\dim(X)$ is even. Suppose that zero does not lie
in the spectrum of $M$. Then the operator 
$\widetilde{d} + \widetilde{d}^*$ is invertible.  (More precisely, it is
invertible as an operator on a Hilbert $C^*_r \Gamma$-module of
differential forms on $M$.) This implies that
$\ind(\widetilde{d} + \widetilde{d}^*)$ vanishes in 
$K_0(C^*_r \Gamma)$.

The higher index theorem says that 
\begin{equation}
\ind(\widetilde{d} + \widetilde{d}^*)
= A (\nu_*([d+d^*])).
\end{equation}
Let $A_\C : K_0(B\Gamma) \otimes \C \rightarrow K_0(C^*_r \Gamma)
\otimes \C$ be
the complexified assembly map. 
Using the isomorphism $K_0(B\Gamma) \otimes \C \cong 
\HH_{even}(B\Gamma; \C)$, 
the higher index theorem implies that in $K_0(C^*_r \Gamma)
\otimes \C$,
\begin{equation}
\ind(\widetilde{d} + \widetilde{d}^*)_\C = A_\C(\nu_*(*L(X))).
\end{equation}
By assumption, $A_\C$ is injective. This gives a contradiction.

Let $T$ be the operator obtained by restricting 
$\widetilde{d} + \widetilde{d}^*$ to
$$\Lambda^{\frac{dim(X)}{2}}(M) \oplus 
\overline{
\widetilde{d}\Lambda^{\frac{dim(X)}{2}}(M)} \oplus
*\overline{
\widetilde{d}\Lambda^{\frac{dim(X)}{2}}(M)}.$$ 
One can show that the other differential forms on $M$ cancel out
when computing the rational
index of $\widetilde{d} + \widetilde{d}^*$, so $T$ 
will have the same index as
$\widetilde{d} + \widetilde{d}^*$. Then the same arguments apply to $T$ to 
give
$0 \in \sigma \left(\triangle_{\frac{dim(X)}{2}}\right)$.

If $\dim(X)$ is odd, consider the even-dimensional manifold
$X^\prime = X \times S^1$ and the group $\Gamma^\prime = \Gamma \times \Z$. 
As the proposition holds for $X^\prime$, it
must also hold for $X$.
\end{pf}

\begin{corollary} \label{SNC}
Let $X$ be a closed Riemannian
manifold. Let $[X] \in \HH_{dim(X)}
(X; \C)$ be its fundamental class. Suppose that there is a surjective 
homomorphism $\pi_1(X) \rightarrow \Gamma$ such that $\Gamma$ satisfies
SNC and the composite map $X \rightarrow B\pi_1(X) \rightarrow B\Gamma$
sends $[X]$ to a nonzero element of
$\HH_{dim(X)}(B\Gamma; \C)$. Let $M$ be the induced normal $\Gamma$-cover
of $X$. Then on $M$,
$0 \in \sigma \left(\triangle_{\frac{dim(X)}{2}}\right)$ 
if $\dim(X)$ is even and 
$0 \in \sigma \left(\triangle_{\frac{dim(X) \pm 1}{2}} \right)$ 
if $\dim(X)$ is odd.
\end{corollary}
\begin{pf}
As the Hirzebruch $L$-class starts out as $L(X) = 1 + \ldots$,
its Poincar\'e dual is of the form $*L(X) = \ldots + [X]$. 
The corollary follows from Proposition \ref{SNC1}.
\end{pf}

\begin{corollary} \label{aspherical}
Let $X$ be a closed aspherical Riemannian
manifold whose fundamental group satisfies SNC.
Then on $\widetilde{X}$,
$0 \in \sigma \left(\triangle_{\frac{dim(X)}{2}}\right)$ 
if $\dim(X)$ is even and 
$0 \in \sigma \left(\triangle_{\frac{dim(X) \pm 1}{2}} \right)$ 
if $\dim(X)$ is odd.
\end{corollary}
\begin{pf}
This follows from Corollary \ref{SNC}.
\end{pf}
\noindent
{\bf Examples : }\\
1. If $X = T^n$ then Corollary \ref{aspherical} is consistent with 
Example 2 of Section \ref{deff}.\\
2. If $X$ is a compact quotient of $H^{2n}$ then
Corollary \ref{aspherical} is consistent with 
Example 3 of Section \ref{deff}.\\
3. If $X$ is a compact quotient of $H^{2n+1}$ then
Corollary \ref{aspherical} is consistent with 
Example 4 of Section \ref{deff}.\\
4. If $X$ is a closed nonpositively-curved locally symmetric space
then Corollary \ref{aspherical} is consistent with the second remark
after Proposition \ref{hyper}.\\

If $X$ is a closed aspherical manifold, it is known that SNC implies
that the rational Pontryagin classes of $X$ are homotopy-invariants
\cite{Kasparov (1988)} and that $X$ does not admit a Riemannian metric
of positive scalar curvature \cite{Rosenberg (1983)}. Thus we see that
these three questions about aspherical manifolds, namely homotopy-invariance 
of rational Pontryagin classes, (non)existence of positive-scalar-curvature
metrics and the zero-in-the-spectrum question,
are roughly all on the same footing. 

If $X$ is a closed aspherical Riemannian manifold, one can ask for which
$p$ one has $0 \in \sigma(\triangle_p)$ on $\widetilde{X}$. 
The case of locally symmetric spaces
is covered by the second remark after Proposition \ref{hyper}.
Another interesting class of aspherical manifolds consists of those with
amenable fundamental group.  By \cite{Cheeger-Gromov (1986)},
$\Ker(\triangle_p) = 0$ for all $p$. By Corollary \ref{Brooks}, 
$0 \in \sigma(\triangle_0)$.
\begin{proposition} \label{nil}
If $X$ is a closed aspherical manifold such that $\pi_1(X)$ has a nilpotent 
subgroup of finite index then
$0 \in \sigma(\triangle_p)$ for all $p \in [0, \dim(X)]$.
\end{proposition}
\begin{pf}
First, $X$ is homotopy-equivalent to an infranilmanifold, that is,
a quotient of the form $\Gamma\backslash G / K$ where $K$ is a finite group,
$G$ is the semidirect product of $K$ and a connected simply-connected
nilpotent Lie group and $\Gamma$ is a discrete cocompact subgroup of $G$ 
\cite[Theorem 6.4]{Farrell-Hsiang (1983)}. We may as well assume that
$X = \Gamma\backslash G / K$. By passing to a finite
cover, we may assume that $K$ is trivial. That is, $X$ is a nilmanifold.  
From \cite[Corollary 7.28]{Raghunathan (1972)}, $\HH^p(X; \C) \cong
\HH^p(g, \C)$, the Lie algebra cohomology of $g$. From
\cite{Dixmier (1955)}, $\HH^p(g, \C) \neq 0$ for all $p \in [0, \dim(X)]$.
Thus for all $p \in [0, \dim(X)]$, $\HH^p(X; \C) \neq 0$.

Now let $\omega$ be a nonzero harmonic $p$-form on $X$. Let $\pi^* \omega$
be its pullback to $\widetilde{X}$. The idea is to construct
low-energy square-integrable $p$-forms on $X$ by multiplying
$\pi^* \omega$ by appropriate functions on $X$. We define the
functions as in \cite[\S 2]{Brooks (1981)}.
Take a smooth
triangulation of $X$ and choose a fundamental domain $F$ for the
lifted triangulation of $\widetilde{X}$. If $E$ is a finite subset of
$\pi_1(X)$, let $\chi_H$ be the characteristic function of
$H = \cup_{g \in E} \: g \cdot F$. Given numbers  $0 < \epsilon_1 <
\epsilon_2 < 1$,
choose a nonincreasing function $\psi \in C^\infty_0([0,\infty))$
which is identically one on $[0, \epsilon_1]$ and identically zero on
$[\epsilon_2, \infty)$. Define a compactly-supported function $f_E$ on 
$\widetilde{X}$ by $f_E(m) = \psi(d(m, H))$.
Then there is a constant $C_1 > 0$, independent of $E$, such that
\begin{equation}
\int_{\widetilde{X}} |df|^2 \le C_1 \: \area(\partial H). 
\end{equation}
Define $\rho_E \in \Lambda^p(\widetilde{X})$ by $\rho_E = f_E \cdot 
\pi^* \omega$. We have
$d \rho_E = df_E \wedge \pi^* \omega$ and $d^* \rho_E = - i(df_E) \: \pi^*
\omega$. As $f_E$ is identically one on $H$, it follows that there is a
constant $C > 0$, independent of $E$, such that
\begin{equation}
\frac{\int_{\widetilde{X}} \left[ |d\rho_E|^2 + |d^* \rho_E|^2 
\right]}{\int_{\widetilde{X}} |\rho_E|^2} \le 
C \: \frac{\area(\partial H)}{\vol(H)}.
\end{equation}
As $\pi_1(X)$ is amenable, by an appropriate choice of $E$
this ratio can be made arbitrarily small.
Thus $0 \in \sigma(\triangle_p)$.
\end{pf}
\noindent
{\bf Question : } Does the conclusion of Proposition \ref{nil} hold if we
only assume that $\pi_1(X)$ is amenable?
\section{Topologically Tame Manifolds}

Another class of manifolds for which one can hope to get some
nontrivial results about the zero-in-the-spectrum question is
given by {\em topologically tame} manifolds, meaning manifolds $M$
which are diffeomorphic to the interior of a compact manifold $N$ with
boundary. If $M$ has finite volume then $\Ker(\triangle_0) \neq 0$, so
we restrict our attention to the infinite volume case. A limited
result is given by Corollary \ref{cor}.

An interesting subclass of topologically tame manifolds consists of those
which are hyperbolic, that is, of constant sectional curvature $-1$.
Complete hyperbolic manifolds are divided into those which are
{\em geometrically
finite} and those which are {\em geometrically infinite}.  Roughly speaking,
$M$ is geometrically finite if its set of ends consists of a finite number of
standard cusps and flares.
\begin{proposition} (Mazzeo-Phillips 
\cite[Theorem 1.11]{Mazzeo-Phillips (1990)})
Let $M$ be an infinite-volume geometrically finite hyperbolic manifold.
If $\dim(M) = 2k$ then $\dim(\Ker(\triangle_k)) = \infty$.
If $\dim(M) = 2k+1$ then $\sigma(\triangle_k) = 
\sigma(\triangle_{k+1}) = [0, \infty)$.
\end{proposition}
The paper \cite{Mazzeo-Phillips (1990)} also computes $\dim(\Ker(\triangle_p))$
for such manifolds.

In general, geometrically infinite hyperbolic manifolds can have
wild end behavior. However, in three dimensions one can show that
the ends have a fairly nice structure. This is used to prove the next
result.

\begin{proposition} (Canary \cite[Theorem A]{Canary (1992)})
If $M$ is a geometrically infinite topologically tame
hyperbolic $3$-manifold then
$0 \in \sigma(\triangle_0)$.
\end{proposition}
\begin{pf}
The method of proof is to show that $M$ is not open at infinity and then
apply Theorem \ref{Buser}.
See \cite{Canary (1992)} for details.
\end{pf}

Thus zero lies in the spectrum of all topologically tame hyperbolic
$3$-manifolds. From Proposition \ref{outside}, the same statement is true for
compactly-supported modifications of such manifolds.

\end{document}